\documentclass[12pt]{amsart}

\usepackage{stmaryrd}

\usepackage{mathabx}
\input bookman.sty
\boldmath

\usepackage{helvet}

\usepackage{graphics}
\usepackage{amssymb}
\usepackage{amsxtra}
\usepackage{amsmath}
\usepackage{mathrsfs}

\usepackage[arrow, matrix]{xy}
\xyoption{frame}

\theoremstyle{plain}

\newtheorem{theo}{Theorem}[section]
\newtheorem{lemm}[theo]{Lemma}
\newtheorem{prop}[theo]{Proposition}
\newtheorem{coro}[theo]{Corollary}
\newtheorem{conj}[theo]{Conjecture}

\theoremstyle{definition}

\newtheorem{defi}[theo]{Definition}
\newtheorem{rema}[theo]{Remark}

\newfont{\rmm}{cmr10 scaled 1000}

\newfont{\itt}{cmsl10 scaled 1000}

\newfont{\rM}{cmr10 scaled 1700}

\newcounter{lemma}[section]

\newcounter{tempcounter}

\newcommand{\lb}{\label}

\newcommand{\rrf}[1]{(\ref{#1})}

\renewcommand{\a}{\alpha}

\newcommand{\g}{\gamma}

\renewcommand{\l}{\lambda}

\renewcommand{\r}{\rho}

\newcommand{\s}{\sigma}

\newcommand{\G}{\Gamma}

\renewcommand{\L}{\Lambda}

\newcommand{\nn}{{\mathbb{N}}}

\newcommand{\qq}{{\mathbb{Q}}}

\newcommand{\zz}{{\mathbb{Z}}}

\newcommand{\TTTT}{{\mathscr{T}}}

\newcommand{\VVVV}{{\mathscr{V}}}

\newcommand{\kkrest}{\begin{picture}(14,14)
\put(00,04){\line(1,0){14}}
\put(00,02){\line(1,0){14}}
\put(06,-4){\line(0,1){14}}
\put(08,-4){\line(0,1){14}}
\end{picture}     }

\newcommand{\krest}{~\kkrest~}

\newcommand{\Ker}{\text{\rm Ker }}

\newcommand{\bere}{\begin{rema}}
\newcommand{\bede}{\begin{defi}}

\renewcommand{\beth}{\begin{theo}}
\newcommand{\bele}{\begin{lemm}}
\newcommand{\bepr}{\begin{prop}}
\newcommand{\beeq}{\begin{equation}}
\newcommand{\bega}{\begin{gather}}
\newcommand{\begaa}{\begin{gather*}}
\newcommand{\been}{\begin{enumerate}}

\newcommand{\bedee}{\begin{defii}}
\newcommand{\bethh}{\begin{theoo}}
\newcommand{\belee}{\begin{lemmm}}
\newcommand{\beprr}{\begin{propp}}

\newcommand{\beco}{\begin{coro}}

\newcommand{\beal}{\begin{aligned}}

\newcommand{\enre}{\end{rema}}

\newcommand{\enco}{\end{coro}}
\newcommand{\enpr}{\end{prop}}
\newcommand{\enth}{\end{theo}}
\newcommand{\enle}{\end{lemm}}
\newcommand{\enen}{\end{enumerate}}
\newcommand{\enga}{\end{gather}}
\newcommand{\engaa}{\end{gather*}}
\newcommand{\eneq}{\end{equation}}
\newcommand{\enal}{\end{aligned}}

\newcommand{\bq}{\begin{equation}}
\newcommand{\bqq}{\begin{equation*}}

\renewcommand{\leq}{\leqslant}
\renewcommand{\geq}{\geqslant}

\newcommand{\wi}{\widetilde}

\newcommand{\ove}{\overline}

\newcommand{\wh}{\widehat}

\newcommand{\sbs}{\subset}

\newcommand{\tens}[1]{\underset{#1}{\otimes}}

\newcommand{\wrt}{with respect to}
\newcommand{\ho}{homomorphism}

\newcommand{\ma}{manifold}

\newcommand{\Prf}{{\it Proof.\quad}}

\newcommand{\chart}{\Phi_p:U_p\to B^n(0,r_p)}
\newcommand{\atlas}{\{\Phi_p:U_p\to B^n(0,r_p)\}_{p\in S(f)}}

\newcommand{\qs}{\hfill\square}

\newcommand{\pa}{\vskip0.1in}

\newcommand{\arrh}[3]
{
\xymatrix{
{#1} \ar[r]^<<<<{#2}  &{#3}
}
}

\newcommand{\arrr}[1]
{\arrh {}{#1}{}}

\newcommand{\arrto}
{\xymatrix{{} \ar@{|-{>}}[r]  & {} } }

\newcommand{\arrinto}
{\xymatrix{{} \ar@{^{(}->}[r]  & {} } }

\newcommand{\llbb}{\llbracket}
\newcommand{\rrbb}{\rrbracket}

\begin{document}

\title[Finiteness of Hofer-Zehnder capacity and loop space homology] 
{Finiteness of $\pi_1$-sensitive Hofer-Zehnder capacity and equivariant loop space homology}
\author{Urs Frauenfelder and  Andrei Pajitnov}
\address{Urs Frauenfelder, 
Universit\"at Augsburg, Institut f\"ur Mathematik,
Universit\"atsstrasse 14, 86159 Augsburg, Germany}
\email{urs.frauenfelder@math.uni-augsburg.de}
\address{Andrei Pajitnov,
Laboratoire Math\'ematiques Jean Leray 
UMR 6629,
Universit\'e de Nantes,
Facult\'e des Sciences,
2, rue de la Houssini\`ere,
44072, Nantes, Cedex}                    
\email{andrei.pajitnov@univ-nantes.fr}
\begin{abstract}
We prove that if $M$ is a closed, connected, oriented, rationally inessential manifold, then
the Hofer-Zehnder capacity of the unit disk bundle of the cotangent bundle of $M$ is finite.
\end{abstract}
\maketitle

\section{Introduction}

In this note we prove the following theorem.
\\ \\
\textbf{Theorem\,A: } \emph{Assume that $M$ is a closed, connected, oriented, rationally inessential manifold. Then
the Hofer-Zehnder capacity of the unit disk bundle of its cotangent bundle is finite}
\\ \\
Actually we prove more. Namely under the assumptions of the theorem we find a \emph{contractible} periodic orbit whose period bounds the Hofer-Zehnder capacity which implies the stronger result  that the $\pi_1$-sensitive Hofer-Zehnder capacity becomes finite. Here we recall that a manifold is called rationally inessential
if the image of its fundamental class under the classifying map to the classifying space of its fundamental group vanishes in homology with rational coefficients. 

To prove our main result we  use a theorem of Irie. This theorem tells us that if the fundamental class of $M$
vanishes in equivariant loop space homology, then the $\pi_1$-sensitive Hofer-Zehnder capacity is finite. What goes on here is that the period of the killer of the fundamental class of $M$ in equivariant symplectic homology
which by a famous result of Viterbo coincides with equivariant loop space homology bounds the $\pi_1$-sensitive Hofer-Zehnder capacity. We then show taking advantage of a deep result of Goodwillie that for
a rationally inessential manifold the fundamental class of $M$ vanishes in equivariant homology. 

Iries theorem holds as well in the nonequivariant set-up. It cannot directly be applied to the nonequivariant
homology of the free loop space, since the fundamental class of $M$ never vanishes in nonequivariant loop space homology. However, it can well happen that it vanishes when one twists the coefficients. This has been examined in \cite{AlbersFrauenfelderOancea}. In particular, it is shown there that if there exists a nontrivial invariant homomorphism from $\pi_2(M)$ to any finite cyclic group, then the $\pi_1$-sensitive Hofer-Zehnder
capacity of the unit disk bundle of its cotangent bundle is finite. Here invariant refers to the action of
the fundamental group of $M$ on $\pi_2(M)$. It is interesting to note that equivariant and nonequivariant methods
lead to quite different topological conditions which guarantee finiteness.

\section{The Hofer-Zehnder capacity}

We briefly recall the definition and general properties of the Hofer-Zehnder capacity. A comprehensive reference is the book by Hofer and Zehnder \cite{hofer-zehnder}.

Assume that $(W,\omega)$ is a symplectic manifold. A smooth function $H \colon W \to \mathbb{R}$ 
gives rise to a vector field $X_H \in \Gamma(TW)$ implicitly defined by
$$dH=\omega(\cdot,X_H)$$
which is referred to as the Hamiltonian vector field of $H$.
\bede\lb{d:simple}
A compactly supported, smooth function $H \colon W \to [0,\infty)$ is called \emph{simple}, if it meets the following conditions
\begin{description}
 \item[(i)] There exists an open nonempty subset $U \subset W$, such that $H|_{U}=\max H$.
 \item[(ii)] The only critical values of $H$ are $0$ and $\max H$.
\end{description}
\end{defi}
Pictorially a simple function looks like a reverted pot. 
\bede\lb{d:admis}
A simple function $H$ is called \emph{admissible}, if its Hamiltonian vector field $X_H$ has no nonconstant orbit of period less than or equal to one. 
\end{defi}
We are now in position to define the Hofer-Zehnder capacity of the symplectic manifold $(W,\omega)$
$$c_{HZ}(W,\omega):=\sup \big\{\max H|\,H \colon W \to [0,\infty)\,\,\textrm{admissible}\big\}.$$
The Hofer-Zehnder capacity has the following properties
\begin{description}
 \item[(i)] Assume there exists a symplectic embedding $(W_1,\omega_1) \hookrightarrow (W_2,\omega_2)$ between symplectic manifolds $(W_1,\omega_1)$ and $(W_2,\omega_2)$, then
$$c_{HZ}(W_1,\omega_1) \leq c_{HZ}(W_2,\omega_2).$$
 \item[(ii)] If $a \neq 0$, then
$$c_{HZ}(W,a\omega)=|a|c_{HZ}(W,\omega).$$
 \item[(iii)] For $r>0$ let
$$B_r:=\{z \in \mathbb{C}^n:|z| < r\}$$
be the $r$-ball, and
$$Z_r:=\{z \in \mathbb{C}^n:|z_1| < r\}$$
be the $r$-cylinder, where $\mathbb{C}^n$ is endowed with its standard symplectic structure, then
$$c_{HZ}(B_r)=r^2 \pi =c_{HZ}(Z_r).$$
\end{description}
Properties (i) and (ii) are straightforward, where property (iii) is highly nontrivial. It implies Gromov's nonsqueezing result \cite{gromov}
\beth\lb{t:nonsqueez}
Assume that $r<r'$, then there exists no symplectic embedding $B_{r'} \hookrightarrow Z_r$. $\qs$
\enth
Another application of the Hofer-Zehnder capacity is Struwes almost sure existence theorem. 
\beth\lb{t:almost}
Assume that $(W,\omega)$ is a compact symplectic manifold with connected boundary such that
$$c_{HZ}(W,\omega)<\infty.$$
If $H \colon W \to \mathbb{R}$ is a smooth function satisfying $H|_{\partial W}>\epsilon$ for some
$\epsilon>0$ such that the interval $(-\epsilon,\epsilon)$ consists of regular values of $H$, then for
almost every $c \in (-\epsilon,\epsilon)$ there exists a periodic orbit on the energy level set
$H^{-1}(c)$.$\qs$
\enth
There is the following variant of the Hofer-Zehnder capacity. Instead of asking just for nonconstant periodic orbits in the definition of admissible Hamiltonian one can ask more generally for contractible periodic orbits. The variant of the Hofer-Zehnder capacity is denoted by $c^0_{HZ}(W,\omega)$ and referred to as the \emph{$\pi_1$-sensitive Hofer-Zehnder capacity}. There is the following obvious inequality 
$$c_{HZ}(W,\omega) \leq c_{HZ}^0(W,\omega).$$
Moreover, one has the following variant of Struwe's theorem, which guarantees not just periodic orbits on almost every energy level but more excitingly contractible periodic orbits, namely
\beth\lb{t:almost2}
Under the assumptions of Theorem~\ref{t:almost}, assume that 
$$c_{HZ}^0(W,\omega)< \infty.$$
Then for almost every $c \in (-\epsilon,\epsilon)$ there exists a contractible periodic orbit on the energy level
$H^{-1}(c)$. $\qs$
\enth
In view of Theorem~\ref{t:almost2} a challenging question is to find conditions which guarantee the finiteness of the $\pi_1$-sensitive Hofer-Zehnder capacity. A first tool to do that is the \emph{displacement energy}.

In order to explain the displacement energy we first need to introduce the Hofer norm of a maybe time dependent Hamiltonian. Suppose that $H \in C_c^\infty(W \times [0,1], \mathbb{R})$ is a compactly supported Hamiltonian depending in addition on a time parameter $t \in [0,1]$. For each $t \in [0,1]$ we then obtain a 
compactly supported function $H_t=H( \cdot,t) \in C^\infty_c(W,\mathbb{R})$. The \emph{Hofer norm} of
$H$ is defined as
$$||H||:=\int_0^1 (\max H_t -\min H_t)dt.$$
If $H$ is time independent this is just the oscillation of $H$. 

A time dependent Hamiltonian gives rise to a time dependent Hamiltonian vector field and we denote by $\phi_H$ the time one map of its flow. 
Suppose that $U \subset W$ is an open subset of a symplectic manifold and therefore itself a symplectic manifold. Then the displacement energy of $U$ in $W$ is defined as
$$e(U):=\inf \{||H||: \phi_H(U)\cap U= \emptyset\}.$$
Note that this quantity depends as well on $W$, i.e.,
$$e(U)=e(U,W).$$
We use the convention that the infimum of the empty set is $\infty$. With this convention $U$ is displaceable in $W$ if and only if $e(U,W)< \infty$. The following Theorem is proved in \cite[Chapter\,5]{hofer-zehnder}.
\beth\lb{t:displace}
Suppose that $U \subset W$ is an open subset and the inclusion homomorphism 
$i_* \colon \pi_1(U) \to \pi_1(W)$ is injective. Then
$$c_{HZ}^0(U) \leq e(U),$$
i.e., if $U$ is displaceable, its $\pi_1$-sensitive Hofer-Zehnder capacity is finite. $\qs$
\enth

\section{Irie's theorem}

We first explain the notion of a Liouville domain. Suppose that $(W, \lambda)$ is a connected exact symplectic manifold, i.e., $\lambda$ is a one-form on $W$ with the property that $\omega=d\lambda$ is a symplectic form. Then one can define a vector field $Y$ on $W$ uniquely determined by $\lambda$ which is implicitly defined by the condition
$$\lambda=\omega(Y, \cdot).$$
Assume that $W$ is compact. Because $W$ admits an exact symplectic form it necesarily has a boundary. 
We say that $(W,\lambda)$ is a Liouville domain if $Y$ points outside at the boundary. This implies that
$(\partial W,\lambda|_{\partial W})$ is a contact manifold. In particular, one can define the Reeb vector field
$R \in \Gamma(T \partial W)$ implicitly by the conditions
$$\lambda(R)=1, \quad d\lambda|_{\partial W}(R,\cdot)=0.$$
To a Liouville domain one can associate its symplectic homology $SH_*(W)$ as well as its
equivariant symplectic homology $SH_*^{S^1}(W)$. Moreover, there are inclusion homomorphisms
$$i_* \colon H_*(W,\partial W) \to SH_*(W), \quad i^{S^1}_* \colon H_*(W,\partial W) \otimes H_*(\mathbb{CP}^\infty) \to SH_*^{S^1}(W).$$
The following theorem is due to Irie \cite{Irie}.
\beth\lb{t:irie}
Suppose that $i_*([W])=0 \in SH_*(W)$ or $i_*^{S^1}([W] \times [pt])=0 \in SH_*^{S^1}(W)$
where $[pt]$ is the class of a point in $H_*(CP^\infty)$, then $c_{HZ}^0(W)<\infty$. $\qs$
\enth
\Prf
For $a>0$ which does not lie in the spectrum of the Reeb flow on $\partial W$ we denote by
$HF^a_*(W)$ the  Floer homology of the Liouville domain associated to a Hamiltonian which grows
with slope $a$ on the completion of the Liouville domain. If $a<b$ the two Floer homologies are related by the inclusion homomorphism
$$HF_*^a(W) \to HF_*^b(W)$$
and the symplectic homology is defined as the direct limit of this directed system. Analogously equivariant symplectic
homology is defined as the direct limit of the equivariant Floer homologies of the Liouville domain for Hamiltonians of
increasing slope on the completion. Moreover, 
the inclusion homomorphism $i_*$ factors through  
$$H_*(W,\partial W) \to HF_*^a(W) \to SH_*(W).$$
 If $i_*([W])=0$, then because the symplectic homology is the direct limit of the Floer homologie $HF^a_*(W)$
there exists $a>0$ such the fundamental class of the Liouville domain already vanishes in $HF_*^a(W)$. However, Irie proved in \cite[Corollary 3.5]{Irie} that this implies that the $\pi_1$-sensitive Hofer-Zehnder
capacity is bounded from above by $a$. 
$\qs$

\section{Essential and inessential manifolds}
\lb{s:ess_and_iness}

In \cite{Gromov1}
M. Gromov introduced the notion
of essential manifold.
To recall the definition we will need some 
terminological conventions.

In  Sections 
\ref{s:ess_and_iness}
--
\ref{s:applic_sympl}
the term {\it manifold } will 
mean {\it closed connected oriented $C^\infty$ manifold}.
For a  connected topological space $X$ we denote 
by $\l_X:X\to B\pi_1(X)$ the 
(homotopy unique) map inducing identity 
isomorphism in $\pi_1$.
Let $R$ be one of the rings $\zz, \ \qq$.

For a  manifold $M$ we denote by 
$[M]_R$ its fundamental class in $H_*(M, R)$ and by 
$\llbb M \rrbb_R$ the element $(\l_X)_*([M]) \in H_*(B\pi_1(X), R)$
(the index $R$ will be omitted if the value of $R$ will be clear from the context).

\bede
A manifold  $M$ is called {\it $R$-essential}
if $\llbb M \rrbb_R \not=0$.
It is called
{\it $R$-inessential}
if $\llbb M \rrbb_R=0$. We abbreviate  {\it $\zz$-inessential}
to {\it inessential};
a $\qq$-inessential manifold is called also {\it rationally inessential}. 
Similar convention will be used for $\zz$-essential and $\qq$-essential 
manifolds. 
\end{defi}

The property of being essential is closely related to the geometry of the manifold.
The main theorem of \cite{Gromov1}
establishes the following inequality
for any essential manifold of dimension $n$:

\begin{equation}\lb{f:syst}
sys_1(M)\leq C_n \big(Vol M\big)^{\frac 1n}
\end{equation}
\noindent
where $sys_1(M)$ is the {\it systolic constant},
that is, the lower bound of the lengths
of closed non contractible curves in $M$.

Later on it turned out that the property of being inessential
is also closely related to the geometry of the underlying manifold,
namely to its macroscopic dimension and to the existence 
of a riemannian metric of positive scalar curvature.

\bede\lb{d:macro-dim}
A manifold $M$ has {\it macroscopic dimension }
less than $k$ if its universal covering $\wi M$ admits a map 
$\phi:\wi M\to X$ to a $k$-dimensional polyhedron $X$
such that the inverse images of points of $x$ are uniformly bounded.
(Here $\wi M$ is endowed with a Riemannian metric, induced from $M$.)
The least number $k$ such that 
$M$
has macroscopic dimension $\leq k$ is called 
{\it macroscopic dimension of} $M$,
and denoted by $\dim_{md} M$.
If $\dim_{md} M=\dim M$
the manifold $M$ is called {\it macroscopically large},
if not, it is called 
{\it macroscopically small}.
\end{defi}
Gromov conjectured \cite{Gromov2} that a compact manifold $M^n$ with 
a Riemannian metric of positive scalar curvature 
satisfies the inequality
$\dim_{md} M \leq \dim M -2.$

This conjecture is still open
(as well as its weaker version 
affirming that a manifold 
with 
a Riemannian metric of positive scalar curvature 
is  macroscopically small).
Inessential manifolds provide a large source of examples
(and counterexamples) to questions related to the Gromov conjecture above.

It is proved in \cite{Dranishnikov1}
that a rationally inessential manifold
$M$ is macroscopically small
if the fundamental group of $M$ is a duality group 
of cohomological dimension not equal to $n+1$.
Dranishnikov \cite{Dranishnikov2} conjectured 
that every rationally inessential manifold $M^n$ 
is macroscopically small. This conjecture was disproved by
M. Marcinkowski, who constructed a rationally inessential 
macroscopically large manifold. The fundamental group 
of this manifold is a finite index subgroup of a Coxeter group.
It would be interesting to have more information about the class of 
groups $G$ for which any rationally inessential manifold $M$
with $\pi_1(M)\approx G$ is 
macroscopically small.

In the two next sections 
we establish a relation between 
the inessentiality property and Hofer-Zehnder capacity for symplectic manifolds.

\section{Some properties of inessential manifolds}
\lb{s:inessential}

Let us begin with some examples of inessential manifolds.

{\bf 1. \ \ }
If $G$ is an $R$-acyclic group
(that is, $H_*(G,R)=0$)
then any manifold $M$ with $\pi_1(M)\approx G$ is obviously 
$R$-inessential. 
The class of $\qq$-acyclic groups is rather large,
it includes for example all Coxeter groups 
(by a theorem of M. Davis, see \cite{MDavis}, p. 302).
There are also many $\zz$-acyclic groups, 
for example the famous Higman group \cite{Higman}
given by the following presentation 
$$
\langle x_i ~|~ x_{i+1}=[x_i, x_{i+1}] \rangle, \ \ i\in \zz/4\zz.
$$
By a theorem of Baumslag, Dyer, and Miller 
\cite{BaumslagDyerMiller}
any group can be embedded into a $\zz$-acyclic group. 

{\bf 2. \ \ }
Let $cd G$ denote the cohomological dimension of $G$.
If $\dim M> cd \pi_1(M)$, then $M$ is clearly inessential.
In particular, this is the case if $\pi_1(M)$ is free, and $\dim M >1$.
Same for manifolds $M$ with $\dim M \geq 3$, if $\pi_1(M)$ is a knot group.

Here are some general properties of inessential manifolds.
(Recall that $R$ stands for $\zz$ or $\qq$.)
\bepr\lb{p:product}
Let $M$ be an $R$-inessential manifold, and $N$ be
any    manifold
of dimension $>1$.
Then $M\times N$ is inessential.
\enpr
\Prf
Let $m=\dim M, n=\dim N$,
put $G=\pi_1(M), H=\pi_1(N)$.
We have 
$$
\llbb M\times N \rrbb = 
\llbb M \rrbb
\otimes 
\llbb N \rrbb
$$
(using the K\"unneth identification of 
$H_m(BG)\tens{R} H_n(BH) $
with a submodule of 
$H_{m+n}(B(G\times H) )$).
The proposition follows. $\qs$
\bepr\lb{p:sum}
Let $M, N$ be 
  manifolds
of same dimension $\geq 3$. Then $M\krest N$ is $R$-inessential
if and only if both $M$ and $N$ are $R$-inessential.
\enpr
\Prf 
Let $m=\dim M$. 
Observe that $B\pi_1(M\krest N)\sim B(G * H) \sim BG\vee BH$
so $H_m(B\pi_1(M\krest N))\approx H_m(BG)\oplus H_M(BH).$
The class $\llbb M\krest N\rrbb$
equals 
$\llbb M\rrbb\oplus
\llbb N\rrbb$
with respect to this decomposition. The proposition follows. $\qs$

\bepr\lb{p:finite-deg}
Let $M, N$ be   manifolds
of same dimension; assume that $M$ is $\qq$-inessential. 
Assume that there is a map $M\to N$ of 
non-zero degree. Then $N$ is $\qq$-inessential.
\enpr
\Prf Obvious. $\qs$

\bepr\lb{p:cov}
Let $p:M\to N$ be a finite covering of manifolds.
Then $N$ is $\qq$-inessential 
if and only if $M$ is $\qq$-inessential.
\enpr 
\Prf The {\it if} part follows from the previous proposition.
To prove the {\it only if} part, 
let $H=\pi_1(M), \ G=\pi_1(N)$. The \ho~ $H\to G$ induced by $P$ is
injective, and we have a commutative diagram,
where both  vertical arrows are finite coverings:

\begin{equation}\lb{f:dia}
 \xymatrix{
M \ar[d]^{p} \ar[r]^{\l_M} & BH \ar[d]^{p'}\\
N \ar[r]^{\l_N} & BG
}
\end{equation}
Recall the transfer \ho~ 
$T:H_*(N)\to H_*(M)$ defined on the chain level by the formula
$T(\s)=\sum \g$ where the sum ranges over all lifts $\g$ of the singular simplex $\s$ to
$M$. The \ho~  $p_*\circ T$ equals identity in $H_*(N)$.
We have also 
$[N]=k\cdot p_*([M])$, where $K$ is the degree of the covering, and  
 $[M]=T([N])$.
The commutativity of the diagram  \rrf{f:dia}
together with the functoriality of the transfer homomorphism
implies that 
$$
\llbb M \rrbb =  T (\llbb N \rrbb ).
$$
The result follows. $\qs$

\section{Applications to symplectic homology theory}
\lb{s:applic_sympl}

Our application of the theory of inessential manifolds to
sympectic homology theory uses Viterbo's 
theorem (see 
\cite{Viterbo} and 
\cite{Viterbo2}, and also 
\cite{AbbondandoloSchwarz} and 
\cite{SalamonWeber}),
which identifies symplectic homology of the cotangent bundle
with homology of the loop space of the underlying manifold.

Let $M$ be a manifold; in the rest of this section all  homology and cohomology groups 
 are with rational coefficients,
 the free loop space $\L M$ is endowed with the canonical
rotation action, and we endow $M$ with the trivial $S^1$-action.
Denote by $\alpha$ the composition
\begin{equation}\lb{f:alfa}
 H_*(M)\arrr {\gamma} H_*^{S^1}(M) \arrr {i_*} H_*^{S^1}(\L M).
\end{equation}

\bepr\lb{p:iness-loops}
If $M$ is $\qq$-inessential, then $\alpha([M])=0$.
\enpr
\Prf
Denote by $\wh\alpha$ the composition of the \ho s 
in cohomology dual to \rrf{f:alfa}.
It suffices to prove that $\wh\a_n:H^n_{S^1}(\L M) \to H^n(M)$ is equal to $0$.
To this end observe that 
$H^*_{S^1}( M)$ is a free module over $H^*(BS^1) \approx \qq[u]$,
therefore the map 
$H^n_{S^1}( M) \to H^n(M)$ is isomorphic to the projection
$$
\bigoplus_{i\in \nn} H^{n-2i}(M) \to H^n(M)
$$
and factors through the inclusion
$$
H^n_{S^1}( M)
\arrinto 
S^{-1} 
H^n_{S^1}( M)
$$
(where the localization is with respect to the multiplicative subset $S=\{u^k~|~ k\in \nn\}$).
Thus it suffices to prove that the following composition vanishes in degree $n$
\begin{equation}\lb{f:compo}
H^*_{S^1}(\L M) \to S^{-1}H^*_{S^1}(\L M) \to S^{-1}H^*_{S^1}(M)
\to H^*(M).
\end{equation}
The sequence \rrf{f:compo}
is obviously functorial in $M$
so the classifying map $\Phi:M\to BG$ induces a commutative diagram

$\xymatrix{
H^*_{S^1}(\L M) \ar[r] & S^{-1}H^*_{S^1}(\L M) \ar[r]&  S^{-1}H^*_{S^1}(M) \ar[r]&  H^*(M) \\
H^*_{S^1}(\L BG) \ar[u] \ar[r] & S^{-1}H^*_{S^1}(\L BG) \ar[u]\ar[r] & S^{-1}H^*_{S^1}(BG) \ar[u]\ar[r]&  H^*(BG) \ar[u]
}$

The second vertical arrow is an isomorphism by T. Goodwillie's theorem
\cite{Goodwillie}, Corollary V.3.3, 
therefore the vanishing of the upper horizontal composition 
follows from vanishing of the homomorphism 
$H^n(BG)\to H^n(M)$. $\qs$

\beco\lb{c:van}
Let $M$ be a $\qq$-inessential manifold. 
Denote by $DM$ and $SM$ the disk bundle (respectively the sphere bundle)
associated with the tangent bundle of $M$.
Then 
$$i_*^{S^1} ([DM]\times [pt])=0\in SH^{S^1}(DM).$$
\enco
\Prf
The following diagram is commutative

\hspace{2cm} $\xymatrix{
H_*(DM, SM) \ar[r] & SH^{S^1}_*(DM) \\
H_*(M) \ar[u]^{\TTTT}_\approx \ar[r] &  H_*^{S^1}(\L M) \ar[u]^{\VVVV}_\approx
}$

(Here $\TTTT$ is the Thom isomorphism, and $\VVVV$ is the Viterbo's isomorphism 
\cite{Viterbo}.) The bottom arrow sends the fundamental class of $M$
to zero by Proposition
\ref{p:iness-loops}
therefore the fundamental class of $DM$ is sent to zero 
by the top arrow. $\qs$ 

Combining Irie's theorem with  Corollary \ref{c:van} 
we deduce a result about the 
finiteness of HZ-capacity for tangent bundles.

\beth\lb{t:iness-capac}
Let $M$ be a $\qq$-inessential manifold.
Then the $\pi_1$-sensitive HZ capacity of $DM$ 
is finite. $\qs$
\enth

In the next Corollary we gathered some immediate 
consequences of Theorem \ref{t:iness-capac}.

\beco\lb{c:list}
Let $M$ be a  manifold.
the $\pi_1$-sensitive HZ capacity of $DM$ 
is finite, if one of the conditions below holds.
\been
\item $\pi_1(M)$ is finite.
\item $\pi_1(M)$ is $\qq$-acyclic.
\item
$\pi_1(M)$ is free and $\dim M >1$.
\item
$M\approx N_1\times N_2$, where $N_1, N_2$ are 
oriented  manifolds, 
$N_1$ is $\qq$-inessential
and $\dim N_2 >0$. $\qs$
\enen 
\enco

\section{On the non-oriented case}
\label{s:non-orient}

Our results can be generalized to  the case of 
non-oriented manifolds.
In this section all manifolds are assumed to be 
closed and connected, but not necessarily oriented. 
For such a \ma~ $M$ 
there is a twisted fundamental class 
$[M]\in H_n(M, w_1)$, where $n=\dim M$ and 
$w_1$ is the local system of groups $\qq$ 
determined by the first Stiefel-Whitney class.
Let $\ove{M}\to M$ be the two-fold  orientation covering of $M$.
The transfer \ho~ sends the twisted homology 
$H_*(M, w_1)$ to the usual homology of $\ove{M}$, and 
we have 
\begin{equation}\lb{f:or-cov}
 T([M])=[\ove{M}].
\end{equation}
\bede\lb{d:iness-non-or}
A non-oriented closed connected manifold $M$ is called 
{\it
$\qq$-inessential}, if the image 
$\llbb M \rrbb \in H_*(BG, w_1)$
of its twisted fundamental class $[M]$ 
equals $0$.
\end{defi}

\bepr\lb{p:nonor-or}
Let $M$ be a non-oriented $\qq$-inessential manifold.
Then $\ove M$ is $\qq$-inessential.
\enpr
\Prf
Let $G=\pi_1(M)$, and $H=\Ker(w_1:G\to \zz/2\zz)$.
We have a commutative diagram 
\begin{equation}\lb{f:diaa}
 \xymatrix{
\ove M \ar[d]^{p} \ar[r]^{\l_{\ove M}} & BH \ar[d]^{p'}\\
M \ar[r]^{\l_M} & BG
}
\end{equation}
and applying the functoriality of the transfer homomorphism
we deduce the proposition. $\qs$

In order to deduce the finiteness of the HZ-capacities for non-oriented essential 
manifolds we will need the following lemma.
\bele\lb{l:HZ-cover}
Let $p:W\to V$ be a finite covering of symplectic manifolds.
Then 
$$
c_{HZ}(W) \geq c_{HZ}(V), \ \ c^0_{HZ}(W) \geq c^0_{HZ}(V).
$$
\enle
\Prf
Let $f$ be an  $HZ$-admissible function on $V$.
Then $f\circ p$ is  an $HZ$-admissible function on $W$,
and the first inequality follows.
Similarly, if $g$ is an  $HZ^0$-admissible function on $V$,
the function $g\circ p$ is  an $HZ$-admissible function on $W$
(indeed, if the Hamiltonian vector field $X_{ g\circ p}$
has a contractible periodic orbit $\g$ of period $\leq 1$
on $W$, then $p\circ\g$ is a contractible periodic orbit for the Hamiltonian 
 vector field $X_{ g}$). The result follows. $\qs$
 
\beco\lb{HZ-nonor}
Let $M$ be a non-oriented $\qq$-inessential manifold.
Then  $c^0_{HZ}(DM)$ is finite. $\qs$
\enco

\section{On  symplectic invariants of unit disc bundles}
\lb{s:unit_bundles}

Let $M$ be a compact closed Riemannian manifold.
The tangent bundle $TM$ has a natural symplectic structure, 
and the symplectic invariants
of the unit disc bundle $DM$ are therefore
invariants of the Riemannian manifold $M$.
The next proposition follows immediately from the definition of
the HZ capacity.

\bepr\lb{p:length}
Let $\g(M)$ denote the minimal length of a non-constant
closed geodesic on $M$.
Then 
$\g(M)\leq C_{HZ}(DM)$.
In particular we have the following lower bound for the HZ capacity of $DM$:
$\r(M)\leq C_{HZ}(M)$
where $\r(M)$ is the injectivity radius of $M$.$\qs$
\enpr

For a non-simply connected manifold $W$
there is a natural refinement of the definition of 
the Hofer-Zehnder capacities (see for example \cite{Irie}).

\bede\lb{d:HZ-refined}
Let $W$ be a symplectic manifold.
Denote by $\langle \pi_1(W)\rangle$ the set of conjugacy classes
of $\pi_1(W)$. 
Let $\Gamma \sbs \langle \pi_1(W)\rangle$
be any non-empty subset of $\pi_1(W)$.
A simple function $H$ is called 
{\it HZ-admissible \wrt~ $\Gamma$}
if the Hamiltonian vector field $X_H$ of $H$ has no non-constant 
periodic orbits $\g$ of period $\leq 1 $ 
with $[\g]\in \Gamma$.
Denote by $C_{HZ}^\Gamma(W)$ 
the upper bound of the set
$$
\{
\max H
~|~
H \ \ {\rm is } \ HZ-{\rm admissible \ \ with \ respect \ to ~ \ } \Gamma
\}
$$
\end{defi}
\noindent
We have obviously 
$$
C_{HZ}^\Gamma(W)
\geq 
C_{HZ}(W)
$$
for any $\Gamma$.
The $\pi_1$-sensitive capacity considered above 
equals by definition 
$C_{HZ}^{\{1\}}(W)$. 

\begin{conj}\lb{conj:iness_gamma}
Let $M$ be an inessential closed compact manifold.
For any $\G\sbs \langle \pi_1(M)\rangle$
the capacity $C_{HZ}^\Gamma(DM)$ is finite.
\end{conj}
\pa

The case of essential manifolds is not covered by the methods of the present paper.
However there is an obvious relation between the $HZ$ capacities and 
systolic constants of $M$.
Namely, 
let $\langle \pi_1(M)\rangle'$
denote the set of all non-trivial 
conjugacy classes of the  fundamental group.
Recall that the  lower bound of lengths of closed non-contractible curves in $M$
is denoted by 
$sys_1(M)$.
We have obviously 
$$
sys_1(M)\leq C_{HZ}\big(DM, \langle \pi_1(M)\rangle'\big).
$$
In view of Gromov's inequality
\rrf{f:syst}
it is natural to suggest the following conjecture

\begin{conj}\lb{conj:syst_capac}
Let $M^n$ be an essential closed Riemannian  manifold
of dimension $n$.
Then
$$
\big(Vol M\big)^{\frac 1n}
\leq D_n
\cdot
C_{HZ}\big(DM, \langle \pi_1(M)\rangle'\big),
$$
where
$D_n$ depend only on $n$.
\end{conj}

\section{Acknowledgements}

The first author would like
to thank the City University of Hong Kong and the Universit\'e de Nantes
for the excellent research atmosphere during the visits March 2015 and November 2015.
The second author is grateful to the City University of Hong Kong
and to the University of Augsburg for the warm hospitality 
during the visits (March 2015 and May 2015) that allowed us
to start and helped to continue our work.

\end{document}